\documentclass[a4paper,twoside,english,DIV12]{scrartcl}
\usepackage[T1]{fontenc}
\usepackage[latin9]{inputenc}
\usepackage{amsthm}
\usepackage{amsmath}
\usepackage{amssymb}
\usepackage{graphicx}
\usepackage{esint}
\usepackage{tikz,tikz-3dplot}
\usepackage{xifthen}									

\makeatletter

 \usepackage{amsthm}
 \theoremstyle{plain}
\newtheorem{thm}{Theorem}[section]
\numberwithin{equation}{section}
\numberwithin{figure}{section}
  \theoremstyle{remark}
  \newtheorem*{rem*}{\protect\remarkname}

\usepackage{enumerate}\usepackage{a4wide}
\usepackage{url}\usepackage{euscript}
\allowdisplaybreaks[1] 
\newenvironment{keywords}{ \noindent\footnotesize\textbf{Keywords and phrases:}}{}

\newenvironment{class}{\noindent\footnotesize\textbf{Mathematics subject classification 2010:}}{}

\usepackage{mathtools}

\usepackage{amsfonts}

\newcommand{\hide}[1]{}

\newcommand*{\restricted}[1]{{\mid_{#1}}}

\newcommand*{\curl}{\operatorname{curl}}

\newcommand*{\ii}{\mathrm{i}}

\DeclareMathAccent{\Circ}{\mathalpha}{operators}{"17}
\newcommand{\interior}[1]{\Circ{#1}}

\renewcommand{\Re}{\operatorname{\mathfrak{Re}}}

\newcommand{\oi}[2]{\left]#1,#2 \right[}

\newcommand{\lci}[2]{\left[#1,#2 \right[}

\renewcommand{\hat}{\widehat}
\renewcommand{\tilde}{\widetilde}
\renewcommand*{\epsilon}{\varepsilon}
\renewcommand*{\rho}{\varrho}

\makeatother

\usepackage{babel}
  \providecommand{\remarkname}{Remark}

\begin{document}

\title{A Note on the Justification of the\\
 Eddy Current Model in Electrodynamics.}

\author{Dirk Pauly%
\thanks{Fakultät für Mathematik, Universität Duisburg-Essen, Campus Essen,
Germany, dirk.pauly@uni-due.de %
} ~ \& Rainer Picard%
\thanks{Institut für Analysis, Technische Universität Dresden, Germany, rainer.picard@tu-dresden.de %
}}
\maketitle
\begin{abstract}
\textbf{Abstract.} The issue of justifying the eddy current approximation
of Maxwell's equations is re-considered in the time-dependent setting.
Convergence of the solution operators is shown in the sense of strong
operator limits. 
\end{abstract}
\begin{keywords}
  electrodynamics, eddy current problem, evolutionary equations, material laws \end{keywords}

\begin{class}
35Q61, 78A25, 78A48
\end{class}\setcounter{section}{-1}

\date{\today}

\section{Introduction}

Somewhat simplified, the topic of this paper is the study of the limit
$\epsilon\to0+$ in the standard Maxwell system in a non-empty open
set $\Omega\subseteq\mathbb{R}^{3}$, which in convenient block operator
matrix notation is given by
\[
\left(\partial_{0}\left(\begin{array}{cc}
\epsilon & 0\\
0 & \mu
\end{array}\right)+\left(\begin{array}{cc}
\sigma & 0\\
0 & 0
\end{array}\right)+\left(\begin{array}{cc}
0 & -\curl\\
\interior\curl & 0
\end{array}\right)\right)\left(\begin{array}{c}
E\\
H
\end{array}\right)=\left(\begin{array}{c}
-J\\
K
\end{array}\right),
\]
with the limit case $\epsilon=0$ being the so-called eddy current
case. Here, following standard physics notation $\partial_{0}$ denotes
the time derivative, $E,$ $H$, denote the electric and magnetic
field, respectively, and $J,\, K$ corresponding external source terms,
$\epsilon$ denotes the dielectricity and $\mu$ the magnetic permeability,
$\sigma$ denotes the conductivity. The overset circle in $\interior\curl$
is supposed to indicate that the so-called electric boundary condition
is imposed on $E$, which makes $\left(\begin{array}{cc}
0 & -\curl\\
\interior\curl & 0
\end{array}\right)$ skew-selfadjoint 
in $L^{2}(\Omega,\mathbb{C}^{3}\times\mathbb{C}^{3})\equiv L^{2}(\Omega,\mathbb{C}^{6})$.
For the purpose of this introduction we may think of $\epsilon,\mu,\sigma$
simply as non-negative real numbers. The approximation question $\epsilon\to0+$
has been considered in the literature commonly in the second order
form, where $H$ has been eliminated from the equations, i.e. the
equations discussed are not Maxwell's equations but rather the abstract
wave equation
\[
\epsilon\partial_{0}^{2}E+\sigma\partial_{0}E+\curl\mu^{-1}\interior\curl\,E=-\partial_{0}J+\curl\mu^{-1}K\eqqcolon-\tilde{J}.
\]
Following common wisdom indeed a time-harmonic regime is assumed,
where $\partial_{0}$ is replaces by $\ii\omega$, where $\omega$
is a real number referred to as frequency, leading to
\[
(\epsilon\omega^{2}-\ii\omega\sigma)E-\curl\mu^{-1}\interior\curl\,E=\tilde{J}.
\]
Due to the selfadjointness of $\curl\mu^{-1}\interior\curl$ in $L^{2}\left(\Omega,\mathbb{C}^{3}\right)$
we have that for $\sigma\not=0$ the number $\epsilon\omega^{2}-\ii\omega\sigma$
is actually in the resolvent set of $\curl\mu^{-1}\interior\curl$
and so the limit $\epsilon\to0+$ is well controlled by the analyticity
of the resolvent 
\[
z\mapsto\big(z-\curl\mu^{-1}\interior\curl\big)^{-1}
\]
on $\mathbb{C}\setminus[0,\infty\,[$~. The situation is less clear
if $\epsilon,\:\sigma$ are allowed to vary -- say they are piece-wise
constant. For example there may be a decomposition of $\Omega$ into
a relative compact, non-empty, open subset $\Omega_{c}\subseteq\Omega$,
where $\sigma=\sigma_{c}>0$ and $\epsilon=\epsilon_{c}\geq0$, and
the rest, where $\sigma=0$ and $\epsilon>0$. For bounded and sufficiently
regular domains $\Omega$ such that a suitable compact embedding result holds,
the limit $\epsilon_{c}\to0+$ can still be established and so a justification
of the eddy current problem can be given. For a survey see \cite{zbMATH05622875}
and the literature quoted there.

A dramatically different situation arises if $\Omega$ is unbounded, e.g.
$\Omega=\mathbb{R}^{3}$.
Then $\sigma=0$ becomes the dominant case with the material behavior
in $\Omega_{c}$ just being a compact perturbation. In this situation
\[
\epsilon\omega^{2}E-\curl\mu^{-1}\interior\curl\,E=\tilde{J}
\]
is our reference case, where now $\epsilon\omega^{2}\in\mathbb{R}\setminus\left\{ 0\right\} $
is always in the continuous spectrum of the operator $\curl\mu^{-1}\interior\curl$.
Thus in contrast to what seems to be claimed in \cite{ammaribuffanedelecjusteddycurmax}
a solution theory in $L^{2}\left(\Omega,\mathbb{C}^{3}\right)$ is
unavailable. The much more demanding issues involved to study such
perturbation problems and to discuss limiting problems is well developed
in connection with the solution theory for exterior boundary value
problems and the study of low-frequency asymptotics in e.g. 
\cite{picardlowfreqmax},
\cite{weckwitschcomlowfreqredwav,weckwitschgenela1,weckwitschgenela2},
\cite{paulydiss,paulyasym}, \cite{pepperldiss},
\cite{ammarinedeleclowfreqemscat}. A comparison between the low-frequency
asymptotics for the full time-harmonic Maxwell's equations and their
eddy current approximation can be found in \cite[Kapitel 5, Satz 5.7]{pepperldiss}.

On the other hand, keeping in mind that time-harmonic problems are
non-physical in so far as they produce infinite energy solutions and
merely serve to describe the time-asymptotic behavior in presence
of a -- perpetual -- time-harmonic forcing, it seems appropriate to
by-pass the above spectral issues altogether by discussing the original
-- physical -- dynamic system directly. This is the perspective of
the following presentation, which is based on concepts derived in
e.g. \cite{PIC_2010:1889,PDE_DeGruyter}. After a brief introduction
into the needed framework we discuss the limit to the eddy current
case in full generality in section \ref{sec:Convergence-to-the}.
In particular, we emphasize that size and boundary regularity of the
underlying domain $\Omega$ play no role in the final result. This is
due to the fact that the classical boundary trace results are superfluous
for the basic solution theory and for obtaining the convergence result.

\section{The Functionalanalytical Framework}

Key to the approach presented here is to consider the closure of differentiation
acting on $C_{1}(\mathbb{R},H)$-functions with compact support, i.e.
functions in $\interior C_{1}(\mathbb{R},H)$, as an operator in $H_{\rho}(\mathbb{R},H)$
with $\rho\in\;]0,\infty[\,$, a weighted $L^{2}$-type space with inner
product 
\[
\langle\varphi\,|\,\psi\rangle_{\rho}\coloneqq\int_{\mathbb{R}}\big\langle\varphi(t)\,|\,\psi(t)\big\rangle_{H}\,\exp(-2\rho t)\, dt,
\]
where $\langle\,\cdot\,|\,\cdot\,\rangle_{H}$ denotes the inner product
of the Hilbert space $H$. The resulting operator 
\[
\partial_{0}:D(\partial_{0})\subseteq H_{\rho}(\mathbb{R},H)\to H_{\rho}(\mathbb{R},H)
\]
turns out, \cite{PIC_2010:1889,PDE_DeGruyter}, to be normal with
\begin{equation}
\Re\partial_{0}=\rho.\label{eq:d0posdef}
\end{equation}
This observation implies that for bounded linear operators
$M_{0}:H\to H$ and $M_{1}:H\to H$, where $M_{0}$ is selfadjoint,
and for a skew-selfadjoint linear operator $A:D(A)\subseteq H\to H$,
which is possibly unbounded, the relation
\[
\Re\big\langle u\,|\,(\partial_{0}M_{0}+M_{1}+A)u\big\rangle_{\rho}=\big\langle u\,|\,(\rho M_{0}+\Re M_{1})u\big\rangle_{\rho}
\]
holds for all $u\in D(\partial_{0})\cap D(A)$. With the assumption that 
\begin{equation}
\rho M_{0}+\Re M_{1}\geq c>0\label{eq:well-pos}
\end{equation}
for all sufficiently large $\rho\in\;]0,\infty[\,$, we obtain that the
closure $\overline{\partial_{0}M_{0}+M_{1}+A}$ and its adjoint $(\partial_{0}M_{0}+M_{1}+A)^{*}=\overline{\partial_{0}M_{0}+M_{1}^{*}-A}$
both have continuous inverses bounded by $1/c$. In particular, the
null spaces of $\overline{\partial_{0}M_{0}+M_{1}+A}$ and $(\partial_{0}M_{0}+M_{1}+A)^{*}$
are both trivial. Thus, we have the following well-posedness result,
see e.g. \cite{PIC_2010:1889,PDE_DeGruyter}.

\begin{thm} \label{thm:well} Let $M_{k}:H\to H$, $k=0,1$, be continuous
linear operators, $M_{0}$ selfadjoint, such that \eqref{eq:well-pos}
holds for some $c\in\;]0,\infty[\,$ and for all $\rho\in\;]\rho_0,\infty[\,$
with $\rho_{0}\in\;]0,\infty[\,$ sufficiently large. Moreover
let $A:D(A)\subseteq H\to H$ be skew-selfadjoint. Then 
\[
(\overline{\partial_{0}M_{0}+M_{1}+A})\, u=f
\]
has for any $f\in H_{\rho}(\mathbb{R},H)$ a unique solution $u\in H_{\rho}(\mathbb{R},H)$.
Furthermore, $u$ depends on $f$ continuously, i.e. 
\[
(\overline{\partial_{0}M_{0}+M_{1}+A})^{-1}:H_{\rho}(\mathbb{R},H)\to H_{\rho}(\mathbb{R},H)
\]
is a continuous linear operator for $\rho\in\;]\rho_0,\infty[\,$. \end{thm}

As a refinement of \eqref{eq:d0posdef} we also find by integration
by parts that for $u\in\interior C_{1}(\mathbb{R},H)$ \big(and so
for $u\in D(\overline{\partial_{0}M_{0}+M_{1}+A})$\big) 
\[
\Re\big\langle u\,|\,\chi_{_{]-\infty,a]}}(\partial_{0}M_{0}+M_{1})u\big\rangle_{\rho}\geq c\big\langle\chi_{_{]-\infty,a]}}u\,|\,\chi_{_{]-\infty,a]}}u\big\rangle_{\rho}.
\]
This yields that we have also causality in the sense of the following
theorem.

\begin{thm}{[}Causality{]} \label{thm:(Causality)} Under the assumptions
of Theorem \ref{thm:well} we have 
\[
\chi_{_{]-\infty,a]}}(\overline{\partial_{0}M_{0}+M_{1}+A})^{-1}=\chi_{_{]-\infty,a]}}(\overline{\partial_{0}M_{0}+M_{1}+A})^{-1}\chi_{_{]-\infty,a]}}
\]
for all sufficiently large $\rho\in\;]0,\infty[\,$ . \end{thm}

We plan to approach the eddy current approximation within this abstract
framework, which simplifies matters in so far as we can deal with
the time-dependent situation under very general assumptions on the
coefficients, which can indeed be operators acting in the underlying
spatial Hilbert space.

\section{\label{sec:Convergence-to-the}Convergence to the Eddy Current Model}

\subsection{Maxwell's Equations with General Material Laws}

Maxwell's equations with a general, simple material law read: 
\begin{align}
\left(\partial_{0}M+N+A\right)\left(\begin{array}{c}
E\\
H
\end{array}\right)=\left(\begin{array}{c}
-J\\
K
\end{array}\right),\qquad A:=\left(\begin{array}{cc}
0 & -\curl\\
\interior\curl & 0
\end{array}\right).\label{eq:Maxwell}
\end{align}
Here $\interior\curl$ is defined as the closure in $L^{2}\left(\Omega,\mathbb{C}^{3}\right)$
of the classical vectoranalytic operation $\curl$ on $C_{1}\left(\Omega,\mathbb{C}^{3}\right)$-vector
fields with compact support in the non-empty open set $\Omega\subseteq\mathbb{R}^{3}$,
which is obviously symmetric in $L^{2}\left(\Omega,\mathbb{C}^{3}\right)$
and therefore indeed closable. We define 
\[
\curl\coloneqq\interior\curl^{*},
\]
which is nothing but the classical weak $L^{2}\left(\Omega,\mathbb{C}^{3}\right)$-curl.
Due to the structure of $A$ as 
\[
A=\left(\begin{array}{cc}
0 & -\interior\curl^{*}\\
\interior\curl & 0
\end{array}\right)
\]
we read off that $A$ is skew-selfadjoint. Since every closed, linear
operator gives rise to a canonical Hilbert space by equipping its
domain with the graph inner product, we have Hilbert spaces 
\[
H(\interior\curl),\,\, H\left(\curl\right)
\]
from the respective domains $D(\interior\curl),\: D\left(\curl\right)$.
One frequently finds already $H(\interior\curl)$ defined
in terms of boundary traces, which unnecessarily limits the applicability
of the results. Even worse, it suggests to the uninitiate or confused
reader that boundary regularity is required to ensure that the physical
model actually works. We note that since
\[
\interior\curl=\curl^{*}\subseteq\curl
\]
we have for $E\in H(\interior\curl)$, with $\left\langle \:\cdot\:|\:\cdot\:\right\rangle _{L^{2}}$
denoting the inner product of $L^{2}\left(\Omega,\mathbb{C}^{3}\right)$,
that 
\[
\langle E\,|\curl\Psi\rangle _{L^{2}}
=\langle \interior\curl\,E\,|\,\Psi\rangle _{L^{2}}
=\langle \curl E\,|\,\Psi\rangle _{L^{2}}
\]
 for all $\Phi\in H\left(\curl\right)$. We read off that conversely
\[
\langle E\,|\curl\Psi\rangle _{L^{2}}
=\langle \curl E\,|\,\Psi\rangle _{L^{2}}\qquad\mbox{ for all }\Psi\in H\left(\curl\right)
\]
characterizes $E\in H(\interior\curl)$. This shows that
\[
E\in H(\interior\curl)
\]
is a suitable generalization of the electric boundary condition for the
topological boundary of arbitrary non-empty open sets.

According to the above abstract framework, the solvability constraint
on the operator coefficients $M=M_{0}$, $N=M_{1}$ is 
\[
\rho M+\Re N\geq c>0
\]
for some real constant $c$ and all sufficiently large $\rho\in\oi0\infty$.
The underlying Hilbert space is $H=L^{2}(\Omega,\mathbb{C}^{6})$.
We recall that causality of the solution operator is also implied
by our general framework.

\subsection{Classical Electrodynamics and the Eddy Current Problem}

On this basis we are now able to discuss the limiting behavior to the
eddy current case. Let 
\[
M_{s}:=\left(\begin{array}{cc}
\epsilon_{s} & 0\\
0 & \mu_{s}
\end{array}\right),\qquad N_{s}:=\left(\begin{array}{cc}
\sigma_{s} & 0\\
0 & 0
\end{array}\right),\qquad s\in[0,1[\,.
\]
Assuming that for some $\hat{\rho}\in\oi0\infty$ we have for all
$\rho\in\oi{\hat{\rho}}\infty$ and all $s\in[0,1[$ 
\[
\rho M_{s}+N_{s}\geq c>0,
\]
we have uniform boundedness for the solution operators in the sense
that 
\[
\big\Vert(\partial_{0}M_{s}+N_{s}+A)^{-1}\big\Vert\leq\frac{1}{c}
\]
for $s\in\lci01$~. On the other hand, we have the following resolvent
equation type result for the solution operators: 
\begin{align*}
 & \qquad(\partial_{0}M_{s}+N_{s}+A)^{-1}-(\partial_{0}M_{0}+N_{0}+A)^{-1}\\
 & =(\partial_{0}M_{s}+N_{s}+A)^{-1}\big((M_{0}-M_{s})\partial_{0}+N_{0}-N_{s}\big)(\partial_{0}M_{0}+N_{0}+A)^{-1}
\end{align*}
If now 
\[
M_{s}\xrightarrow{s\to0+}M_{0},\quad N_{s}\xrightarrow{s\to0+}N_{0}\qquad\text{strongly in }
L^{2}(\Omega,\mathbb{C}^{6}),
\]
we read off that we have 
\begin{align*}
(\partial_{0}M_{s}+N_{s}+A)^{-1}F\xrightarrow{s\to0+}(\partial_{0}M_{0}+N_{0}+A)^{-1}F\qquad\mbox{for every }
F\in D(\partial_{0}).
\end{align*}
Due to the uniform boundedness of the solution operators, however,
we can use the density of $D(\partial_{0})$ in $H_{\rho}\left(\mathbb{R},L^{2}(\Omega,\mathbb{C}^{6})\right)$
and the above $H_{\rho}\left(\mathbb{R},L^{2}(\Omega,\mathbb{C}^{6})\right)$-convergence
for elements in $D(\partial_{0})$. In fact we get 
\[
(\partial_{0}M_{s}+N_{s}+A)^{-1}F\xrightarrow{s\to0+}(\partial_{0}M_{0}+N_{0}+A)^{-1}F\qquad\mbox{for every }
F\in H_{\rho}\left(\mathbb{R},L^{2}(\Omega,\mathbb{C}^{6})\right)
\]
by the principle of uniform boundedness%
\footnote{Indeed, more explicitly, for $F\in H_{\rho}\big(\mathbb{R},L^{2}(\Omega,\mathbb{C}^{6})\big)$
and $\tilde{F}\in D(\partial_{0})$ we see 
\begin{align*}
 & \qquad\big|(\partial_{0}M_{s}+N_{s}+A)^{-1}F-(\partial_{0}M_{0}+N_{0}+A)^{-1}F\big|_{\rho}\\
 & \leq\big|(\partial_{0}M_{s}+N_{s}+A)^{-1}\tilde{F}-(\partial_{0}M_{0}+N_{0}+A)^{-1}\tilde{F}\big|_{\rho}\\
 & \qquad+\big|(\partial_{0}M_{s}+N_{s}+A)^{-1}(F-\tilde{F})\big|_{\rho}+\big|(\partial_{0}M_{0}+N_{0}+A)^{-1}(\tilde{F}-F)\big|_{\rho}\\
 & \leq\big|(\partial_{0}M_{s}+N_{s}+A)^{-1}\tilde{F}-(\partial_{0}M_{0}+N_{0}+A)^{-1}\tilde{F}\big|_{\rho}+\frac{2}{c}|F-\tilde{F}|_{\rho},
\end{align*}
from which the desired convergence result follows by first choosing
$\tilde{F}\in D(\partial_{0})$ to make the last term sufficiently
small (independently of $s\in\lci01$~) and then, for this fixed
choice of $\tilde{F}$, we choose $s_{0}\in\oi01$ sufficiently small
to make the first term sufficiently small for all $s\in\oi0{s_{0}}$\,.%
}, i.e. strong convergence of the solution operators. In cases, where
$\epsilon_{0}$ vanishes, i.e. 
$\epsilon_{0}\restricted{L^{2}(\Omega)^{3}}=0$, we have the case
of the eddy current approximation. We note that the usually considered
case assumes $\Omega=\mathbb{R}^{3}$.

%

Let us conclude with some remarks:
\begin{itemize}
\item 
The above rationale clearly also works for completely general sequences
$(M_{s})_{s}$, $(N_{s})_{s}$ with $M_{s}^{*}=M_{s}$ of continuous linear operators
converging strongly to $M_{0}$ and $N_{0}$ with 
\[
\rho M_{s}+\Re N_{s}\geq c>0
\]
for some $c\in\;]0,\infty[\,$ and all $s\in\lci01$~. We have focused
on the classical eddy current context to make the approach more tangible.
\item 
Again we emphasize that our results depend 
\begin{itemize}
\item
neither on the size (bounded or unbounded)
\item
nor on the topology (genus, Betti-numbers)
\item
nor on the regularity (no regularity is assumed) 
\end{itemize}
of the underlying domain resp. non-empty open set $\Omega$.
Moreover, our methods extend immediately to domains 
resp. non-empty open sets $\Omega\subset\mathbb{R}^{N}$,
$N\in\mathbb{N}$, or even to Riemannian manifolds $\Omega$
by replacing the $\curl$-operators by the exterior resp. co-derivative.
\item 
Our results remain valid even if mixed boundary conditions are considered.
We just have to modify the skew-selfadjoint unbounded linear operator $A$ by
\[
A=\left(\begin{array}{cc}
0 & -\interior\curl_{\Gamma_{1}}^{*}\\
\interior\curl_{\Gamma_{1}} & 0
\end{array}\right)
=\left(\begin{array}{cc}
0 & -\interior\curl_{\Gamma_{2}}\\
\interior\curl_{\Gamma_{1}} & 0
\end{array}\right).
\]
Here the boundary $\Gamma:=\partial\,\Omega$ is decomposed into, let's say,
two relative open disjoint subsets $\Gamma_{1}\not=\Gamma$ 
and $\Gamma_{2}:=\Gamma\setminus\overline{\Gamma_{1}}$.
Following our definitions from above, we define $\interior\curl_{\Gamma_{1}}$ 
as the closure in $L^{2}\left(\Omega,\mathbb{C}^{3}\right)$
of the $\curl$-operator acting on the restrictions to $\Omega$
of $C_{1}\left(\mathbb{R}^3,\mathbb{C}^{3}\right)$-vector fields 
having compact support in $\mathbb{R}^3$ bounded away from the boundary part $\Gamma_{1}$
as well as
\[
\interior\curl_{\Gamma_{2}}\coloneqq\interior\curl_{\Gamma_{1}}^{*}.
\]
Once more, the structure of $A$ shows that $A$ is skew-selfadjoint. 
\item 
It is also clear that for uniform convergence of the coefficients
we get uniform convergence of the solution operators in the sense
that
\[
\partial_{0}^{-1}\left(\partial_{0}M_{s}+N_{s}+A\right)^{-1}\xrightarrow{s\to0+}\partial_{0}^{-1}\left(\partial_{0}M_{0}+N_{0}+A\right)^{-1}
\]
in $\mathcal{L}\left(H_{\rho}(\mathbb{R},L^{2}(\Omega,\mathbb{C}^{6})),H_{\rho}(\mathbb{R},L^{2}(\Omega,\mathbb{C}^{6}))\right).$
\end{itemize}


\subsection{A Realistic Case}

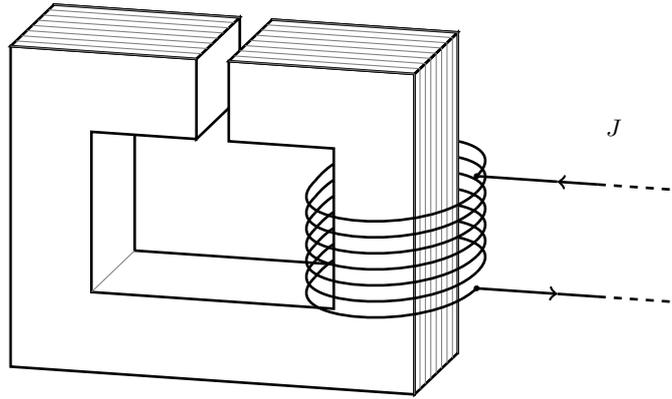
\begin{figure}
\centering
\tdplotsetmaincoords{75}{105}
\begin{tikzpicture}[scale=1.1,tdplot_main_coords]
	\def\gdist{0.2} 
	\def\fdist{-2}  
	\def\sdist{0.2} 
	\foreach \z in {0,\sdist,...,1.2}{
		\draw [line width=1pt] (\fdist/2,1,\z-\sdist-0.7) .. controls (\fdist-1.2,1,\z-0.7) and (\fdist-1.2,3,\z-0.7) .. (\fdist/2,3,\z-0.7);
	}
	\draw [line width=1pt] (\fdist,-2.5,-2) -- (\fdist,2.5,-2) -- (\fdist,2.5,2) -- (\fdist,\gdist,2) -- (\fdist,\gdist,1) -- (\fdist,1.5,1) -- (\fdist,1.5,-1) -- (\fdist,-1.5,-1) -- (\fdist,-1.5,1) -- (\fdist,-\gdist,1) -- (\fdist,-\gdist,2)-- (\fdist,-2.5,2) -- cycle;
	\fill [white] (0,-\gdist,1) -- (0,-\gdist,2) -- (\fdist,-\gdist,2) -- (\fdist,-\gdist,1) -- cycle;
	\draw [line width=1pt] (0,-\gdist,1) -- (0,-\gdist,2) -- (\fdist,-\gdist,2) -- (\fdist,-\gdist,1) -- cycle;
	\fill [white,opacity=1] (0,-2.5,-2) -- (0,2.5,-2) -- (0,2.5,2) -- (0,\gdist,2) -- (0,\gdist,1) -- (0,1.5,1) -- (0,1.5,-1) -- (0,-1.5,-1) -- (0,-1.5,1) -- (0,-\gdist,1) -- (0,-\gdist,2)-- (0,-2.5,2) -- cycle;
	\draw [line width=1pt] (0,-2.5,-2) -- (0,2.5,-2) -- (0,2.5,2) -- (0,\gdist,2) -- (0,\gdist,1) -- (0,1.5,1) -- (0,1.5,-1) -- (0,-1.5,-1) -- (0,-1.5,1) -- (0,-\gdist,1) -- (0,-\gdist,2)-- (0,-2.5,2) -- cycle;
	\fill [white] (0,2.5,-2) -- (\fdist,2.5,-2) -- (\fdist,2.5,2) -- (0,2.5,2) -- cycle;
	\fill [white] (0,2.5,2) -- (\fdist,2.5,2) -- (\fdist,\gdist,2) -- (0,\gdist,2) -- cycle;
	\fill [white] (0,-2.5,2) -- (\fdist,-2.5,2) -- (\fdist,-\gdist,2) -- (0,-\gdist,2) -- cycle;
	\draw [line width=1pt] (0,2.5,-2) -- (\fdist,2.5,-2) -- (\fdist,2.5,2) -- (0,2.5,2) -- cycle;
	\draw [line width=1pt] (0,2.5,2) -- (\fdist,2.5,2) -- (\fdist,\gdist,2) -- (0,\gdist,2) -- cycle;
	\draw [line width=1pt] (0,-2.5,2) -- (\fdist,-2.5,2) -- (\fdist,-\gdist,2) -- (0,-\gdist,2) -- cycle;
	\foreach \z in {0,-0.25,...,\fdist}{
		\draw [gray] (\z,2.5,-2) -- (\z,2.5,2);
		\draw [gray] (\z,2.5,2) -- (\z,\gdist,2);
		\draw [gray] (\z,-2.5,2) -- (\z,-\gdist,2);
	}
	\draw [gray] (0,-1.5,-1) -- (\fdist,-1.5,-1);
	\foreach \z in {-\sdist,0,\sdist,...,1}{
		\draw [line width=1pt] (\fdist/2,1,\z-0.7) .. controls (1.2,1,\z-0.7) and (1.2,3,\z-0.7) .. (\fdist/2,3,\z-0.7);
	}
	\fill (\fdist/2,3,0.5) circle (1pt); 
	\draw [line width=1pt] (\fdist/2,3,0.5) -- (\fdist/2,4,0.5);
	\draw [line width=1pt,<-] (\fdist/2,4,0.5) -- (\fdist/2,4.5,0.5);
	\draw [line width=1pt,dashed] (\fdist/2,4.5,0.5) -- (\fdist/2,5.5,0.5);
	\fill (\fdist/2,3,-\sdist-0.7) circle (1pt); 
	\draw [line width=1pt,->] (\fdist/2,3,-\sdist-0.7) -- (\fdist/2,4,-\sdist-0.7);
	\draw [line width=1pt] (\fdist/2,4,-\sdist-0.7) -- (\fdist/2,4.5,-\sdist-0.7);
	\draw [line width=1pt,dashed] (\fdist/2,4.5,-\sdist-0.7) -- (\fdist/2,5.5,-\sdist-0.7);
	\draw (\fdist/2,4.7,1) node[anchor=south]{\footnotesize$J$};
\end{tikzpicture}
\caption{Laminated iron core in air.}
\label{manu}
\end{figure}

For illustrational purposes we conclude with a realistic example,
where the above limit situation occurs, e.g. an electromagnetic field in the
presence of a laminated iron core in air (possibly with an air gap). Whereas
in air the standard Maxwell equations are used, in the iron core the
eddy current model is frequently assumed, see Figure \ref{manu}.
A possible, simple description would be that $\mu>0$ is a constant and
$\epsilon$ and $\sigma$ are piece-wise constant with
\[
\epsilon=\begin{cases}
\epsilon_{\textsf{air}} & \mbox{ in air},\\
\epsilon_{\textsf{lam}} & \mbox{ in the insulating parts of the laminated iron core,}\\
\epsilon_{\textsf{cor}} & \mbox{ in the metal parts of the laminated iron core,}
\end{cases}
\]
and
\[
\sigma=\begin{cases}
0 & \mbox{ in air},\\
0 & \mbox{ in the insulating parts of the laminated iron core,}\\
\sigma_{\textsf{cor}} & \mbox{ in the metal parts of the laminated iron core,}
\end{cases}
\]
for $\epsilon_{\textsf{air}},\:\epsilon_{\textsf{lam}},\:\epsilon_{\textsf{cor}},\:\sigma_{\textsf{cor}}$ positive numbers. In the above we have established that replacing the 
-- relative to $\sigma_{\textsf{cor}}$ -- 
small value of the dielectricity $\epsilon_{\textsf{cor}}$ can
indeed be replaced by zero, i.e. $\epsilon_{\textsf{cor}}=0$. 
In this situation the approximation result holds for any $\rho\in\oi0\infty$~. 

\vspace*{8mm}
\textsf{\Large\bfseries Acknowledgement}\\
We cordially thank Immanuel Anjam for creating the picture.

\end{document}